\documentclass[12pt,bezier]{article}
\usepackage{times}
\usepackage{booktabs}
\usepackage{pifont}
\usepackage{floatrow}
\usepackage{caption}
\usepackage{mathrsfs}
\usepackage{amsmath}
\usepackage{amsfonts,amsthm,amssymb,mathrsfs,bbding}
\usepackage{txfonts}
\usepackage{graphics,multicol}
\usepackage{graphicx}
\usepackage{color}
\usepackage{amssymb}
\allowdisplaybreaks[3]
\usepackage{caption}
\usepackage{cite}
\usepackage{latexsym,bm}
\usepackage{indentfirst}
\usepackage{color}
\usepackage[colorlinks=true,anchorcolor=blue,filecolor=blue,linkcolor=blue,urlcolor=blue,citecolor=blue]{hyperref}
\usepackage{extarrows}
\usepackage{cite}
\usepackage{latexsym,bm}
\usepackage{mathtools}
\usepackage{enumerate}
\evensidemargin=\oddsidemargin\topmargin=-1.5cm

\newtheorem{thm}{Theorem}[section]
\newtheorem{prob}{Problem}[section]
\newtheorem{claim}{Claim}

\newtheorem{lem}{Lemma}[section]
\newtheorem{cor}{Corollary}[section]

\theoremstyle{definition}

\addtocounter{section}{0}

\begin{document}
\title{Fractional matching, factors and spectral radius in graphs involving minimum degree\footnote{Supported by National Natural Science Foundation of China
{(Nos. 11971445 and 12171440)} and Natural Science Foundation of Henan Province (No. 202300410377).}}
\author{{\bf Jing Lou$^{a}$}, {\bf Ruifang Liu$^{a}$}\thanks{Corresponding author.
E-mail addresses: rfliu@zzu.edu.cn (R. Liu), loujing\_23@163.com (J. Lou), aoguoyan@163.com (G. Ao)},
{\bf Guoyan Ao$^{a, b}$}\\
{\footnotesize $^a$ School of Mathematics and Statistics, Zhengzhou University, Zhengzhou, Henan 450001, China} \\
{\footnotesize $^b$ School of Mathematics and Statistics, Hulunbuir University, Hailar, Inner Mongolia 021008, China}}
\date{}

\date{}
\maketitle
{\flushleft\large\bf Abstract}
A fractional matching of a graph $G$ is a function $f:E(G)\rightarrow [0, 1]$ such that for any $v\in V(G)$, $\sum_{e\in E_{G}(v)}f(e)\leq1$, where
$E_{G}(v)=\{e\in E(G): e~ \mbox{is incident with} ~v~\mbox{in}~G\}$.
The fractional matching number of $G$ is $\mu_{f}(G)=\mathrm{max}\{\sum_{e\in E(G)}f(e):f$ is a fractional matching of $G\}$.
Let $k\in (0,n)$ is an integer. In this paper, we prove a tight lower bound of the spectral radius to guarantee $\mu_{f}(G)>\frac{n-k}{2}$ in a graph with minimum degree $\delta,$ which implies the result on the fractional perfect matching due to Fan et al. [Discrete Math. 345 (2022) 112892].

For a set $\{A, B, C, \ldots\}$ of graphs, an $\{A, B, C, \ldots\}$-factor of a graph $G$ is defined to be a spanning subgraph of $G$ each component of which is isomorphic to one of $\{A, B, C, \ldots\}$.
We present a tight sufficient condition in terms of the spectral radius for the existence of a $\{K_2, \{C_k\}\}$-factor in a graph with minimum degree $\delta,$ where $k\geq 3$ is an integer. Moreover, we also provide a tight spectral radius condition for the existence of a $\{K_{1, 1}, K_{1, 2}, \ldots , K_{1, k}\}$-factor with $k\geq2$ in a graph with minimum degree $\delta,$ which generalizes the result of Miao et al. [Discrete Appl. Math. 326 (2023) 17-32].
\begin{flushleft}
\textbf{Keywords:} Fractional matching, Factor, Spectral radius, Minimum degree

\end{flushleft}
\textbf{AMS Classification:} 05C50; 05C35

\section{Introduction}
Let $G$ be a simple and undirected graph with vertex set $V(G)$ and edge set $E(G).$
The order and size of $G$ are denoted by $|V(G)|=n$ and $|E(G)|=e(G)$, respectively.
We denote by $\delta(G)$ and $\Delta(G)$ the minimum degree and the maximum degree of $G,$ respectively.
Let $i(G)$ be the number of isolated vertices of $G$.
For a vertex $v\in V(G)$, define $E_{G}(v)=\{e\in E(G): e~ \mbox{is incident with} ~v~\mbox{in}~G\}$.
For a vertex subset $S$ of $G$, we denote by $G-S$ and $G[S]$ the subgraph of $G$ obtained from $G$ by deleting the vertices in $S$ together with their incident edges and the subgraph of $G$ induced by $S$, respectively.
For an edge subset $M$ of $G$, let $G[M]$ be the subgraph of $G$ induced by $M$.
We use $C_n, K_n$ and $K_{1, n-1}$ to denote the cycle of order $n$, the complete graph of order $n$ and the star with $n$ vertices, respectively.
Let $G_1$ and $G_2$ be two vertex-disjoint graphs.
We denote by $G_{1}+G_{2}$ the disjoint union of $G_1$ and $G_2$.
The join $G_1\vee G_2$ is the graph obtained from $G_{1}+G_{2}$ by adding all possible edges between $V(G_1)$ and $V(G_2)$.
For undefined terms and notions in this paper, one can refer to \cite{Bondy2008}.

For a graph $G$, we use $A(G)$ to denote the adjacency matrix of $G.$
The largest eigenvalue of $A(G)$, denoted by $\rho(G)$, is called the {\it spectral radius} of $G$.
It is well known that $A(G)$ is nonnegative and irreducible for a connected graph $G$. According to the Perron-Frobenius theorem,
there exists a unique positive unit eigenvector corresponding to $\rho(G)$, which is called the {\it Perron vector} of $G.$

A {\it fractional matching} of $G$ is a function $f:E(G)\rightarrow [0, 1]$ such that for any $v\in V(G)$, $\sum_{e\in E_{G}(v)}f(e)\leq1$.
The {\it fractional matching number} of $G$ is $\mu_{f}(G)=\mathrm{max}\{\sum_{e\in E(G)}f(e):f$ is a fractional matching of $G\}$.
A fractional matching of $G$ is called a {\it fractional perfect matching} if $\sum_{e\in E(G)}f(e)=\frac{n}{2}.$

\begin{thm}[Scheinerman and Ullman \cite{Scheinerman1997}]\label{thm1}
Let $G$ be a graph of order $n.$ Each of the following holds.
\begin{enumerate}[(i)]
\item Any fractional matching satisfies $\mu_{f}(G)\leq \frac{n}{2}.$
\item $2\mu_{f}(G)$ is an integer.
\end{enumerate}
\end{thm}

Let $i(G-S)$ be the number of isolated vertices of graph $G-S$.
Recall the following result due to Scheinerman and Ullman \cite{Scheinerman1997}.

\begin{thm}[Scheinerman and Ullman \cite{Scheinerman1997}]\label{thm2} (The fractional Berge-Tutte formula)
Let $G$ be a graph of order $n$. Then
$$\mu_{f}(G)=\frac{1}{2}(n-\mathrm{max}\{i(G-S)-|S|:\forall S\subseteq V(G)\}).$$
\end{thm}

Inspired by the work of Scheinerman and Ullman \cite{Scheinerman1997}, we first consider the following spectral extremal problem in this paper.

\begin{prob}\label{p1}
Let $G$ be a graph of order $n$ with minimum degree $\delta$ and $\mu_{f}(G)\leq\frac{n-k}{2}$,
where $k\in (0,n)$ is an integer.
What is the maximum spectral radius of $G?$ Moreover, characterize the corresponding spectral extremal graphs.
\end{prob}

Focusing on Problem \ref{p1}, we prove the first main result.

\begin{thm}\label{main1}
Let $G$ be a graph of order $n\geq \mathrm{max}\{6\delta+5k+1, 5\delta+k^2+4k+1\}$ with minimum degree $\delta$,
where $k\in (0,n)$ is an integer. If
$$\rho(G)\geq \rho(K_{\delta} \vee (K_{n-2\delta-k}+(\delta+k)K_1)),$$
then $\mu_{f}(G)> \frac{n-k}{2}$ unless $G\cong K_{\delta} \vee (K_{n-2\delta-k}+(\delta+k)K_1).$
\end{thm}

Taking $k=1$ in Theorem \ref{main1} and Combining Theorem \ref{thm1}, one can obtain a tight sufficient condition based on the spectral radius for the existence of a fractional perfect matching.

\begin{cor}[Fan, Lin and Lu \cite{Fan2022}]
Let $G$ be a graph of order $n\geq 6\delta+6$ with minimum degree $\delta$.
If $$\rho(G)\geq \rho(K_{\delta}\vee(K_{n-2\delta-1}+(\delta+1)K_{1})),$$
then $G$ contains a fractional perfect matching unless $G\cong K_{\delta}\vee(K_{n-2\delta-1}+(\delta+1)K_{1}).$
\end{cor}

Amahashi and Kano \cite{Amahashi1982} introduced a new notion concerning factors. For a set $\{A, B, C, \ldots\}$ of graphs,
an $\{A, B, C, \ldots\}$-factor of a graph $G$ is a spanning subgraph of $G$ each component of which is isomorphic to one of $\{A, B, C, \ldots\}$.
That is, if $F$ is an $\{A, B, C, \ldots\}$-factor of graph $G$,
then $F$ is a subgraph of $G$ such that $V(F)=V(G)$ and each component of $F$ is contained in $\{A, B, C, \ldots\}$.

In the past few years, many researchers focused on finding conditions for the existence of a perfect matching.
It is easy to see that $G$ contains a perfect matching if and only if $G$ has a $\{K_2\}$-factor.
Tutte \cite{Tutte1947} provided a sufficient and necessary condition for a graph to have a $\{K_2\}$-factor.
Later, Tutte \cite{Tutte1952} presented the following sufficient and necessary condition for the existence of a $\{K_2, \{C_k\}\}$-factor.

\begin{thm}[Tutte \cite{Tutte1952}]\label{thm3}
Let $G$ be a graph and $k\geq 3$ be an integer.
Then $G$ has a $\{K_2, \{C_k\}\}$-factor if and only if
$i(G-S)\leq |S|$
for every $S\subseteq V(G)$.
\end{thm}

\begin{prob}
Let $G$ be a graph of order $n$ with minimum degree $\delta$ which has no a $\{K_2, \{C_k\}\}$-factor.
What is the maximum spectral radius of $G?$ Moreover, characterize all the extremal graphs.
\end{prob}

Using Theorem \ref{thm3}, we prove a tight spectral condition to guarantee the existence of a $\{K_2, \{C_k\}\}$-factor in a graph with minimum degree $\delta$.

\begin{thm}\label{main2}
Let $G$ be a graph of order $n\geq 5\delta +6$ with minimum degree $\delta$ and $k\geq 3$ is an integer. If
$$\rho(G)\geq \rho(K_{\delta}\vee (K_{n-2\delta -1}+(\delta +1)K_1)),$$
then $G$ has a $\{K_2, \{C_k\}\}$-factor unless $G\cong (K_{\delta}\vee (K_{n-2\delta -1}+(\delta +1)K_1))$.
\end{thm}

Akiyama, Avis and Era \cite{Akiyama1980} presented a sufficient and necessary condition for a graph to contain a $\{K_{1, 1}, K_{1,2}\}$-factor.
Soon after, Amahashi and Kano \cite{Amahashi1982} generalized the previous result and posed the following sufficient and necessary condition for a graph to contain a $\{K_{1, 1}, K_{1,2}, \ldots, K_{1,k}\}$-factor,
where $k\geq 2$ is an integer.

\begin{thm}[Amahashi and Kano \cite{Amahashi1982}]\label{thm4}
Let $G$ be a graph and $k\geq 2$ be an integer.
Then $G$ has a $\{K_{1, 1}, K_{1, 2}, \ldots , K_{1, k}\}$-factor if and only if
$i(G-S)\leq k|S|$
for every $S\subseteq V(G)$.
\end{thm}

\begin{prob}
Let $G$ be a graph of order $n$ with minimum degree $\delta$ which has no a $\{K_{1, 1}, K_{1, 2}, \ldots , K_{1, k}\}$-factor,
where $k\geq 2$ is an integer. What is the maximum spectral radius of $G?$
Moreover, characterize spectral extremal graphs.
\end{prob}

\begin{figure}[H]
\centering
{\includegraphics[width=0.32\textwidth]{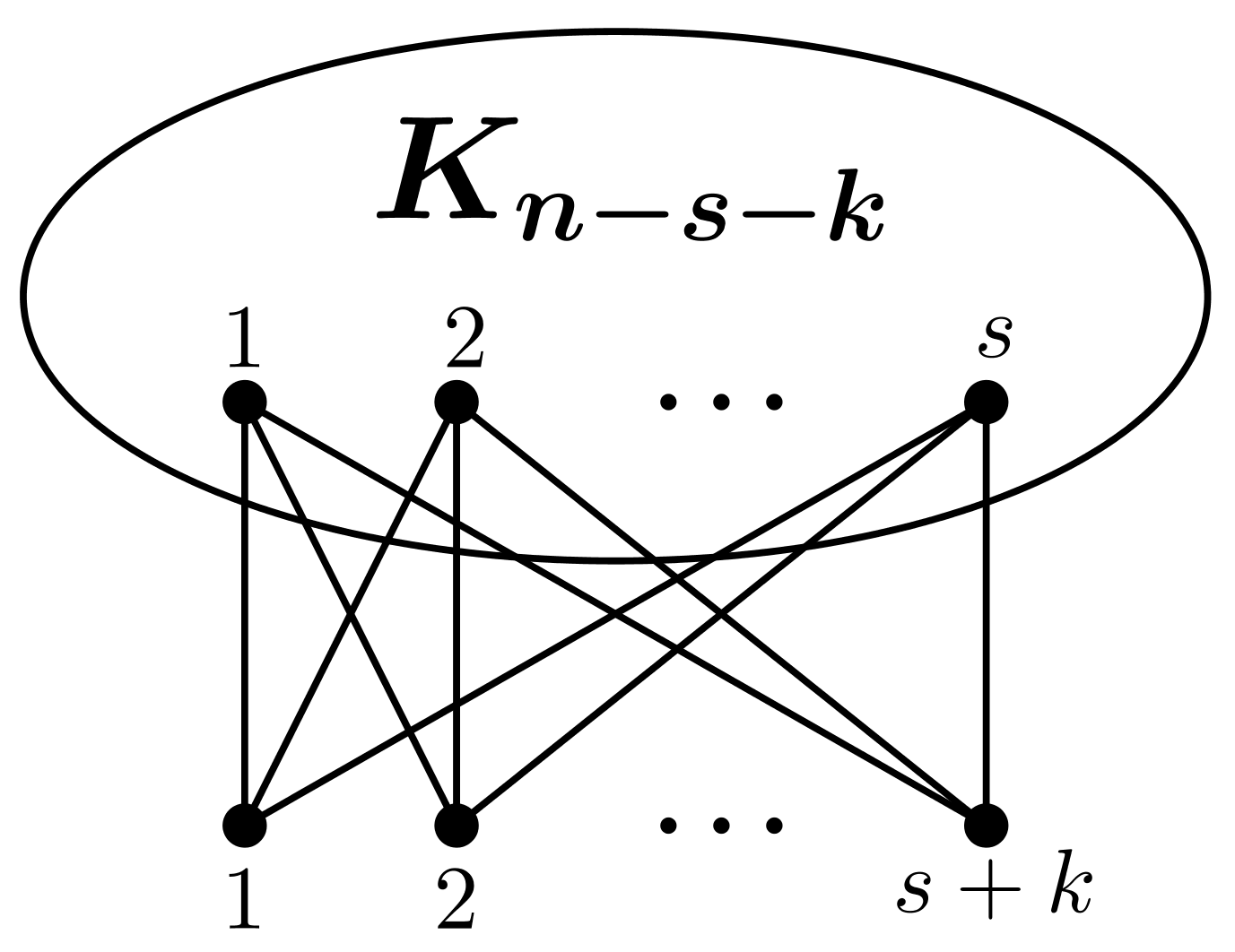}}\\
\caption{The graph $K_{s}\vee (K_{n-2s-k}+(s+k)K_{1}).$}
\label{fig1}
\end{figure}

Miao and Li \cite{Miao2023} established a lower bound on the size (resp. the spectral radius) of
$G$ to guarantee that $G$ contains a $\{K_{1, 1}, K_{1, 2}, \ldots , K_{1, k}\}$-factor. Next we generalize their result to
graphs with minimum degree $\delta$ and prove the following result.

\begin{thm}\label{main3}
Let $G$ be a graph of order $n\geq 3(k+1)\delta +5$ and minimum degree $\delta$, where $k\geq 2$ is an integer. If
$$\rho(G)\geq \rho(K_{\delta}\vee (K_{n-k\delta-\delta-1}+(k\delta+1)K_1)),$$
then $G$ has a $\{K_{1, 1}, K_{1, 2}, \ldots, K_{1, k}\}$-factor unless $G\cong (K_{\delta}\vee (K_{n-k\delta-\delta-1}+(k\delta+1)K_1)$.
\end{thm}

\section{Proof of Theorem \ref{main1}}

Before proving our main result, we first present a preliminary result which will be used in our subsequent arguments.

Let $A=(a_{ij})$ and $B=(b_{ij})$ be two $n\times n$ matrices.
Define $A\leq B$ if $a_{ij}\leq b_{ij}$ for all $i$ and $j$, and define $A< B$ if $A\leq B$ and $A\neq B$.

\begin{lem}[Berman and Plemmons \cite{Berman1979}, Horn and Johnson \cite{Horn1986}]\label{le2}
Let $A=(a_{ij})$ and $B=(b_{ij})$ be two $n\times n$ matrices with the spectral radius $\lambda(A)$ and $\lambda(B)$.
If $0\leq A\leq B$, then $\lambda(A)\leq \lambda(B)$.
Furthermore, if $B$ is irreducible and $0\leq A < B$, then $\lambda(A)<\lambda(B)$.
\end{lem}

Now, we shall give the proof of Theorem \ref{main1}.

\medskip
\noindent  \textbf{Proof of Theorem \ref{main1}.}
Let $G$ be a graph of order $n\geq \mbox{max}\{6\delta+5k+1, 5\delta+k^2+4k+1\}$ and minimum degree $\delta.$
Assume to the contrary that $\mu_{f}(G)\leq \frac{n-k}{2}.$ By Theorem \ref{thm2}, there exists a vertex subset $S\subseteq V(G)$ such that $i(G-S)-|S|\geq k.$
Let $T$ be the set of isolated vertices in $G-S$, $|S|=s$ and $|T|=t$. Then $t=i(G-S)\geq s+k.$ Note that $s+t\leq n.$ Then
$s\leq\frac{n-k}{2}$. Since $N_G(T) \subseteq S$, we have $s\geq \delta.$ It is obvious that $G$ is a spanning subgraph of $G'=K_{s}\vee (K_{n-2s-k}+(s+k)K_1)$~(see Fig. \ref{fig1}). According to Lemma \ref{le2}, we have
\begin{equation}\label{equ1}
\rho(G)\leq \rho(G'),
\end{equation}
where equality holds if and only if $G\cong G'$. Note that $\delta\leq s\leq\frac{n-k}{2}.$ Next we shall divide the proof into two cases according to the value of $s$.

\vspace{1.8mm}
\noindent\textbf{Case 1.} $s=\delta$.
\vspace{1mm}

Then $G'\cong K_{\delta} \vee (K_{n-2\delta-k}+(\delta+k)K_{1})$.
By (\ref{equ1}), we conclude that
$$\rho(G)\leq \rho(K_{\delta} \vee (K_{n-2\delta-k}+(\delta+k)K_{1})).$$
By the assumption $\rho(G)\geq \rho(K_{\delta}\vee (K_{n-2\delta-k}+(\delta+k)K_{1})),$ then we have
$G\cong K_{\delta}\vee (K_{n-2\delta-k}+(\delta+k)K_{1}))$.
Note that
$$\mathrm{max}\{i(G-S)-|S|:\forall S\subseteq V(G)\}=i(G-V(K_{\delta}))-|V(K_{\delta})|=k.$$
By Theorem \ref{thm2}, $\mu_{f}(G)=\frac{n-k}{2}.$
Hence $G\cong K_{\delta}\vee (K_{n-2\delta-k}+(\delta+k)K_{1}))$.

\vspace{1.8mm}
\noindent\textbf{Case 2.} $\delta+1\leq s\leq\frac{n-k}{2}.$
\vspace{1mm}

Recall that $G'=K_{s}\vee (K_{n-2s-k}+(s+k)K_1)$. We can partition the vertex set of $G'$ as $V(G')=V(K_{s})\cup V(K_{n-2s-k})\cup V((s+k)K_{1}).$
Let $V(K_{s})=\{u_1, u_2, \ldots, u_s\}, V(K_{n-2s-k})=\{v_1, v_2, \ldots, v_{n-2s-k}\}$ and $V((s+k)K_{1})=\{w_1, w_2, \ldots, w_{s+k}\}.$
Then
$$A(G')=\left[
\begin{array}{ccc}
(J-I)_{s\times s}&J_{s\times (n-2s-k)}&J_{s\times (s+k)}\\
J_{(n-2s-k)\times s}&(J-I)_{(n-2s-k)\times (n-2s-k)}&O_{(n-2s-k)\times(s+k)}\\
J_{(s+k)\times s}&O_{(s+k)\times (n-2s-k)}&O_{(s+k)\times (s+k)}
\end{array}
\right],
$$
where $J$ denotes the all-one matrix, $I$ denotes the identity square matrix, and $O$ denotes the zero matrix.
Let $x$ be the Perron vector of $A(G')$, and let $\rho'=\rho(G')$.
By symmetry, $x$ takes the same value (say $x_1$, $x_2$ and $x_3$) on the vertices of $V(K_{s}), V(K_{n-2s-k})$ and $V((s+k)K_{1})$, respectively.
According to $A(G')x=\rho'x$, we have
$$\rho'x_3=sx_1.$$
Note that $\rho'>0$. Then
\begin{equation}\label{equ5}
x_3=\frac{sx_1}{\rho'}.
\end{equation}

Define $G^*=K_{\delta}\vee (K_{n-2\delta-k}+(\delta+k)K_{1})$.
Suppose that $E_1=\{v_{i}w_{j}|1\leq i\leq n-2s-k, \delta+k+1\leq j\leq s+k\}\cup \{w_{i}w_{j}|\delta+k+1\leq i\leq s+k-1, i+1\leq j\leq s+k\}$ and $E_2=\{u_{i}w_{j}|\delta+1\leq i\leq s, 1\leq j\leq \delta+k\}$. It is easy to see that $G^*\cong G'+E_1-E_2.$
Note that $A(G^*)$ can be obtained by replacing $s$ with $\delta$ in $A(G')$.
Similarly, let $y$ be the Perron vector of $A(G^*)$ and $\rho^*=\rho(G^*)$.
By symmetry, $y$ takes the same value on the vertices of $V(K_{\delta})$, $V(K_{n-2\delta-k})$ and $V((\delta+k)K_1),$ respectively.
We denote the entry of $y$ by $y_1$, $y_2$ and $y_3$ corresponding to the vertices in the above three vertex subsets, respectively.
By $A(G^*)y=\rho^*y$, we have
\begin{gather}
\rho^*y_2=\delta y_1+(n-2\delta-k-1)y_2,\label{equ6}\\
\rho^*y_3=\delta y_1.\label{equ7}
\end{gather}
Note that $G^*$ contains $K_{n-2\delta-k}$ as a proper subgraph.
By Lemma \ref{le2}, we have $\rho^*>\rho(K_{n-2\delta-k})=n-2\delta-k-1$.
Putting (\ref{equ7}) into (\ref{equ6}) and combining $\rho^*>n-2\delta-k-1$, we have
\begin{equation}\label{equ8}
y_2=\frac{\rho^*y_3}{\rho^*-(n-2\delta-k-1)}.
\end{equation}
Clearly, $G'=K_{s}\vee (K_{n-2s-k}+(s+k)K_{1})$ is not a complete graph, and hence $\rho'<n-1.$
Note that $\delta+1 \leq s \leq \frac{n-k}{2}$.
\begin{claim}\label{claim1}
$\rho^*>\rho'.$
\end{claim}

\medskip
\noindent  \textbf{Proof.}
Suppose to the contrary that $\rho'\geq \rho^*$. Then
\begin{eqnarray*}
&&y^{T}(\rho^*-\rho')x\\
&=&y^{T}(A(G^*)-A(G'))x\\
&=& \sum_{uv\in E_{1}}(x_{u}y_{v}+x_{v}y_{u})-\sum_{uv\in E_{2}}(x_{u}y_{v}+x_{v}y_{u})\\
&=&\sum_{i=1}^{n-2s-k}\sum_{j=\delta+k+1}^{s+k}(x_{v_{i}}y_{w_j}+x_{w_j}y_{v_i})+
\sum_{i=\delta+k+1}^{s+k-1}\sum_{j=i+1}^{s+k}(x_{w_i}y_{w_j}+x_{w_j}y_{w_i})
-\sum_{i=\delta+1}^{s}\sum_{j=1}^{\delta+k}(x_{u_i}y_{w_j}+x_{w_j}y_{u_i})\\
&=&(n-2s-k)(s-\delta)(x_2y_2+x_3y_2)+(s-\delta-1)(s-\delta)x_3y_2-(s-\delta)(\delta+k)(x_1y_3+x_3y_2)\\
&=&(s-\delta)[(n-2s-k)(x_2y_2+x_3y_2)+(s-\delta-1)x_3y_2-(\delta+k)(x_1y_3+x_3y_2)]\\
&=&(s-\delta)[(n-s-2\delta-2k-1)x_3y_2+(n-2s-k)x_2y_2-(\delta+k)x_1y_3]\\
&\geq&(s-\delta)[(n-s-2\delta-2k-1)x_3y_2-(\delta+k)x_1y_3],
\end{eqnarray*}
where the last inequality follows from the fact that $n\geq 2s+k$, $x>0$ and $y>0$.
By (\ref{equ5}) and (\ref{equ8}), we have
\begin{eqnarray*}
&&y^{T}(\rho^*-\rho')x\\
&\geq&(s-\delta)x_1y_3\left[\frac{s\rho^*(n-s-2\delta-2k-1)}{\rho'(\rho^*-(n-2\delta-k-1))}-(\delta+k)\right]\\
&=&\frac{(s-\delta) x_1y_3}{\rho'(\rho^*-(n-2\delta-k-1))}[s\rho^*(n-s-2\delta-2k-1)-\rho'(\delta+k)(\rho^*-(n-2\delta-k-1))]\\
&=&\frac{\rho^*(s-\delta) x_1y_3}{\rho'(\rho^*-(n-2\delta-k-1))}[s(n-s-2\delta-2k-1)+(\delta+k)(\frac{\rho'}{\rho^*}(n-2\delta-k-1)-\rho')].
\end{eqnarray*}
Since $\rho^*\leq \rho'<n-1$, then
\begin{eqnarray*}
&&y^{T}(\rho^*-\rho')x\\
&\geq&\frac{\rho^*(s-\delta) x_1y_3}{\rho'(\rho^*-(n-2\delta-k-1))}[s(n-s-2\delta-2k-1)+(\delta+k)(n-2\delta-k-1-\rho')]\\
&>&\frac{\rho^*(s-\delta) x_1y_3}{\rho'(\rho^*-(n-2\delta-k-1))}[s(n-s-2\delta-2k-1)-(\delta+k)(2\delta+k)]\\
&=&\frac{\rho^*(s-\delta) x_1y_3}{\rho'(\rho^*-(n-2\delta-k-1))}[-s^2+(n-2\delta-2k-1)s-2\delta^2-3k\delta-k^2]\\
&\triangleq& \frac{\rho^*(s-\delta) x_1y_3}{\rho'(\rho^*-(n-2\delta-k-1))}f(s).
\end{eqnarray*}
Note that $f(s)$ is a convex function on $s$.
We claim that $f(s)>0$ for $\delta+1 \leq s \leq \frac{n-k}{2}.$

In fact, it suffices to prove that $f(\delta+1)>0$ and $f(\frac{n-k}{2})>0.$
Note that $n\geq5\delta+k^2+4k+1$, $\delta\geq0$ and $k\geq 1$. By a direct calculation, we have
\begin{eqnarray*}
f(\delta+1)&=&(\delta+1)n-5\delta^2-5k\delta -5\delta-k^2-2k-2\\
&\geq&(\delta+1)(5\delta+k^2+4k+1)-5\delta^2-5k\delta -5\delta-k^2-2k-2\\
&=& (k^{2}-k+1)\delta+2k-1\\
&\geq&2k-1> 0,
\end{eqnarray*}
and
\begin{eqnarray*}
f(\frac{n-k}{2})&=&\frac{n^2}{4}-(\delta+k+\frac{1}{2})n-2\delta^2-2k\delta-\frac{k^2}{4}+\frac{k}{2}\\
&\triangleq&g(n).
\end{eqnarray*}
Note that the symmetry axis of $g(n)$ is $x=2\delta+2k+1<6\delta+5k+1$.
It follows that $g(n)$ is a monotonically increasing function on $n\in[6\delta+5k+1,+\infty)$.
Note that $\delta\geq 0$ and $k\geq1$. Then
\begin{eqnarray*}
f(\frac{n-k}{2})=g(n)&\geq& g(6\delta+5k+1)\\
&=&\delta^2+(2k-1)\delta +k^2-\frac{k}{2}-\frac{1}{4}\\
&\geq&k^2-\frac{k}{2}-\frac{1}{4}>0.
\end{eqnarray*}
This implies that $f(s)>0$ for $\delta+1\leq s\leq \frac{n-k}{2}.$

Combining $s\geq\delta+1$ and $\rho^*>n-2\delta-k-1,$ we obtain that $y^{T}(\rho^*-\rho')x>0$.
Hence $\rho^*>\rho',$ which contradicts the assumption $\rho'\geq\rho^*$. \hspace*{\fill}$\Box$

By Claim \ref{claim1} and (\ref{equ1}), we have
$$\rho(G)\leq\rho(G')<\rho(G^*),$$
a contradiction. This completes the proof. \hspace*{\fill}$\Box$

\section{Proof of Theorem \ref{main2}}

Let $M$ be the following $n\times n$ matrix
\[
M=\left(\begin{array}{ccccccc}
M_{1,1}&M_{1,2}&\cdots &M_{1,m}\\
M_{2,1}&M_{2,2}&\cdots &M_{2,m}\\
\vdots& \vdots& \ddots& \vdots\\
M_{m,1}&M_{m,2}&\cdots &M_{m,m}\\
\end{array}\right),
\]
whose rows and columns are partitioned into subsets $X_{1}, X_{2},\ldots ,X_{m}$ of $\{1,2,\ldots, n\}$.
The quotient matrix $R(M)$ of the matrix $M$ (with respect to the given partition)
is the $m\times m$ matrix whose entries are the
average row sums of the blocks $M_{i,j}$ of $M$.
The above partition is called {\it equitable}
if each block $M_{i,j}$ of $M$ has constant row (and column) sum.

\begin{lem}[Brouwer and Haemers \cite{Brouwer2011}, Godsil and Royle \cite{Godsil2001}, Haemers \cite{Haemers1995}]\label{le3}
Let $M$ be a real symmetric matrix and let $R(M)$ be its equitable quotient matrix.
Then the eigenvalues of the quotient matrix $R(M)$ are eigenvalues of $M$.
Furthermore, if $M$ is nonnegative and irreducible, then the spectral radius of the quotient matrix $R(M)$ equals to the spectral
radius of $M$.
\end{lem}

Now we are in a position to present the proof of Theorem \ref{main2}.

\begin{figure}
\centering
\includegraphics[width=0.32\textwidth]{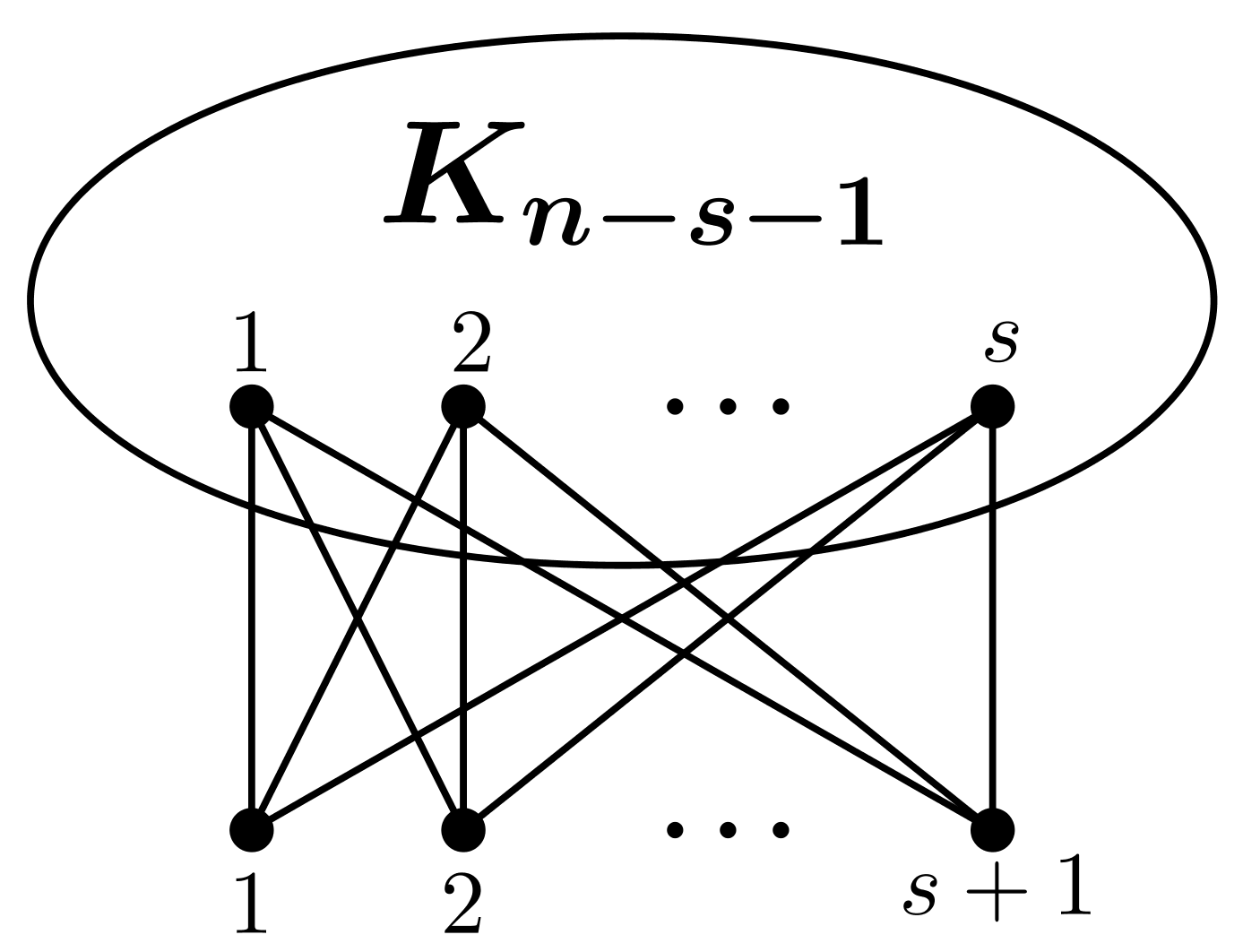}\\
\caption{The graph $K_{s}\vee (K_{n-2s-1}+(s+1)K_{1}).$
}\label{fig2}
\end{figure}

\medskip
\noindent  \textbf{Proof of Theorem \ref{main2}.}
Assume that $G$ has no a $\{K_2, \{C_k\}\}$-factor for $n\geq 5\delta +6$ and $k\geq 3$.
By Theorem \ref{thm3}, there exists a vertex subset $S$ of $V(G)$ such that $i(G-S)\geq |S|+1$.
Let $|S|=s$.
Then $G$ is a spanning subgraph of $G'=K_s\vee (K_{n-2s-1}+(s+1)K_{1})$ (see Fig. \ref{fig2}).
Note that $G$ has the minimum degree $\delta.$ Then $s\geq \delta$. By Lemma \ref{le2}, we have
\begin{equation}\label{equ11}
\rho(G)\leq \rho(G'),
\end{equation}
with equality holding if and only if $G\cong G'$. Then we will divide the proof into the following two cases.

\vspace{1.5mm}
\noindent\textbf{Case 1.} $s=\delta$.
\vspace{1mm}

Then $G'\cong K_{\delta} \vee (K_{n-2\delta-1}+(\delta+1)K_{1})$.
By (\ref{equ11}), we have
$$\rho(G)\leq \rho(K_{\delta} \vee (K_{n-2\delta-1}+(\delta+1)K_{1})).$$
Since $\rho(G)\geq \rho(K_{\delta} \vee (K_{n-2\delta-1}+(\delta+1)K_{1}))$, then
$G\cong K_{\delta} \vee (K_{n-2\delta-1}+(\delta+1)K_{1}).$
Note that the vertices of $(\delta+1)K_1$ are only adjacent to $\delta$ vertices of $K_{n-\delta-1}.$
Then $K_{\delta} \vee (K_{n-2\delta-1}+(\delta+1)K_{1})$ contains no a spanning subgraph each component of which is contained in $\{K_2, \{C_k\}\}$.
This implies that $K_{\delta} \vee (K_{n-2\delta-1}+(\delta+1)K_{1})$ has no a $\{K_2, \{ C_k\}\}$-factor.
Hence $G\cong K_{\delta} \vee (K_{n-2\delta-1}+(\delta+1)K_{1})$.

\vspace{1.5mm}
\noindent\textbf{Case 2.} $s\geq \delta+1.$
\vspace{1mm}

Let $R_{s}$ be an equitable quotient matrix of the adjacency matrix $A(G')$ with respect to the partition
$(V(K_{s}), V(K_{n-2s-1}), V((s+1)K_1))$.
One can see that
$$R_{s}=\left[
\begin{array}{ccc}
s-1&n-2s-1&s+1\\
s&n-2s-2&0\\
s&0&0
\end{array}
\right].
$$
By a simple calculation, the characteristic polynomial of $R_{s}$ is
$$P(R_{s},x)=x^3-(n-s-3)x^2-(s^2+n-2)x-2s^3+ns^2-4s^2+ns-2s.$$
By Lemma \ref{le3}, we know that $\rho(G')=\lambda_{1}(R_{s})$ is the largest root of the equation $P(R_{s},x)= 0$.
Note that $A(K_{\delta}\vee (K_{n-2\delta-1}+(\delta+1)K_{1}))$ has the equitable quotient matrix $R_{\delta}$,
which is obtained by replacing $s$ with $\delta$ in $R_{s}$.
Similarly, by Lemma \ref{le3}, $\rho(K_{\delta}\vee (K_{n-2\delta-1}+(\delta+1)K_{1}))=\lambda_{1}(R_{\delta})$
is the largest root of $P(R_{\delta},x)=0$.
Then
$$P(R_{s},x)-P(R_{\delta},x)=(s-\delta)(x^2-(\delta+s)x-2s^2+ns-2\delta s-4s-2\delta^{2}+n\delta -4\delta+n-2).$$
We claim that $P(R_{s},x)-P(R_{\delta},x)>0$ for $x\in [n-\delta-2, +\infty).$
In fact, let $$f(x)=x^2-(\delta+s)x-2s^2+ns-2\delta s-4s-2\delta^{2}+n\delta -4\delta+n-2.$$
Note that $s\geq \delta+1.$ Then we only need to prove $f(x)>0$ on $x\in [n-\delta-2,+\infty)$.
Note that $n\geq 2s+1$. Then $\delta+1\leq s\leq \frac{n-1}{2}$ and
the symmetry axis of $f(x)$ is $x=\frac{\delta+s}{2}\leq s \leq n-\delta-2.$
This implies that $f(x)$ is increasing with respect to $x\in [n-\delta-2, +\infty)$.
Since $s\leq \frac{n-1}{2}$, $n\geq 5\delta+6$ and $\delta \geq 0$, then
\begin{eqnarray*}
f(x)&\geq& f(n-\delta-2)\\
&=&-2s^2-(\delta+2)s+n^2-2\delta n-3n+2\delta+2\\
&\geq&-2(\frac{n-1}{2})^2-(\delta+2)(\frac{n-1}{2})+n^2-2\delta n-3n+2\delta+2\\
&=& \frac{1}{2}n^2-(\frac{5\delta}{2}+3)n+\frac{5\delta}{2}+\frac{5}{2}\\
&\geq& \frac{1}{2}(5\delta+6)^2-(\frac{5\delta}{2}+3)(5\delta+6)+\frac{5\delta}{2}+\frac{5}{2}\\
&=& \frac{5}{2}\delta+\frac{5}{2}> 0.
\end{eqnarray*}
It follows that $P(R_{s},x)>P(R_{\delta},x)$ for $x\geq n-\delta-2.$
Note that $K_{\delta} \vee (K_{n-2\delta-1}+(\delta+1)K_1)$ contains $K_{n-\delta-1}$ as a proper subgraph.
By Lemma \ref{le2}, we have
$$\rho(K_{\delta} \vee (K_{n-2\delta-1}+(\delta+1)K_1))> \rho(K_{n-\delta-1})=n-\delta-2,$$
and so $\lambda_{1}(R_{s})< \lambda_{1}(R_{\delta}).$
This means that $\rho(G')<\rho(K_{\delta}\vee (K_{n-2\delta-1}+(\delta+1)K_{1})).$
Therefore,
$$\rho(G) \leq \rho(G') <\rho(K_{\delta}\vee (K_{n-2\delta-1}+(\delta+1)K_{1})),$$
a contradiction. The proof is completed.
\hspace*{\fill}$\Box$

\section{Proof of Theorem \ref{main3}}

\begin{figure}
\centering
\includegraphics[width=0.32\textwidth]{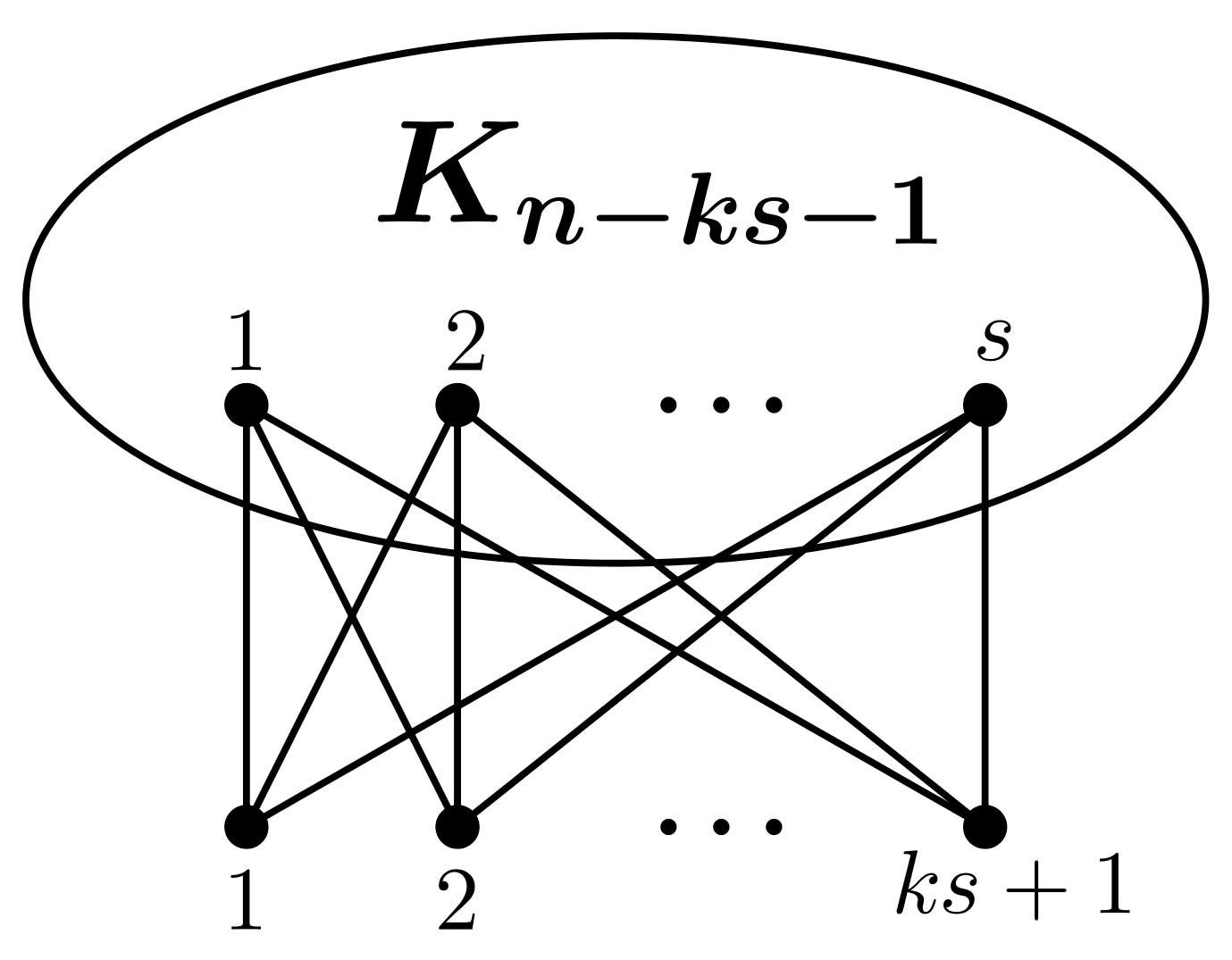}\\
\caption{The graph $K_{s}\vee(K_{n-ks-s-1}+(ks+1)K_{1})$.
}\label{fig3}
\end{figure}

\medskip
\noindent  \textbf{Proof of Theorem \ref{main3}.}
Assume that $G$ is a graph of order $n$ and minimum degree $\delta$, which contains no a $\{K_{1,1}, K_{1,2},\ldots, K_{1,k}\}$-factor for $k\geq 2$ and $n\geq 3(k+1)\delta+5$.
By Theorem \ref{thm4}, there exists a vertex subset $S$ of $V(G)$ such that $i(G-S)\geq k|S|+1$.
Let $|S|=s$. Then $G$ is a spanning subgraph of $G'=K_{s}\vee (K_{n-ks-s-1}+(ks+1)K_{1})$ (see Fig. \ref{fig3}).
Obviously, $s\geq \delta$. By Lemma \ref{le2}, we have
\begin{equation}\label{for1}
\rho(G)\leq \rho(G'),
\end{equation}
where equality holds if and only if $G\cong G'$.
Next we shall divide the proof into the following two cases.

\vspace{1.5mm}
\noindent\textbf{Case 1.} $s=\delta$.
\vspace{1mm}

Then $G'\cong K_{\delta}\vee (K_{n-k\delta-\delta-1}+(k\delta+1)K_{1})$. According to (\ref{for1}), we conclude that
$$\rho(G)\leq \rho(K_{\delta}\vee (K_{n-k\delta-\delta-1}+(k\delta+1)K_{1})).$$
By the assumption,
we have $G\cong K_{\delta}\vee (K_{n-k\delta-\delta-1}+(k\delta+1)K_{1})$.
Note that the vertices of $(k\delta+1)K_{1}$ are only adjacent to $\delta$ vertices of $K_{n-k\delta-1}$.
Then $K_{\delta}\vee (K_{n-k\delta-\delta-1}+(k\delta+1)K_{1})$ has no a spanning subgraph each component of which is contained in
$\{K_{1,1}, K_{1,2},\ldots, K_{1,k}\}$, which implies that
$K_{\delta}\vee (K_{n-k\delta-\delta-1}+(k\delta+1)K_{1})$ has no a $\{K_{1,1}, K_{1,2},\ldots, K_{1,k}\}$-factor.
Hence $G\cong K_{\delta}\vee (K_{n-k\delta-\delta-1}+(k\delta+1)K_{1})$.

\vspace{1.5mm}
\noindent\textbf{Case 2.}  $s\geq \delta+1.$
\vspace{1mm}

Note that $G'=K_{s}\vee (K_{n-ks-s-1}+(ks+1)K_{1})$.
We divide $V(G')$ into three parts: $V(K_{s})$, $V(K_{n-ks-s-1})$ and $V((ks+1)K_{1})$.
Let $V(K_{s})=\{u_1, u_2, \ldots, u_s\}$, $V(K_{n-ks-s-1})=\{v_1, v_2, \ldots, v_{n-ks-s-1}\}$ and
$V((ks+1)K_{1})=\{w_1, w_2, \ldots, w_{ks+1}\}.$
Then
$$A(G')=\left[
\begin{array}{ccc}
(J-I)_{s\times s}&J_{s\times (n-ks-s-1)}&J_{s\times (ks+1)}\\
J_{(n-ks-s-1)\times s}&(J-I)_{(n-ks-s-1)\times (n-ks-s-1)}&O_{(n-ks-s-1)\times(ks+1)}\\
J_{(ks+1)\times s}&O_{(ks+1)\times (n-ks-s-1)}&O_{(ks+1)\times (ks+1)}
\end{array}
\right],
$$
where $J$ denotes the all-one matrix, $I$ denotes the identity square matrix, and $O$ denotes the zero matrix.
Let $x$ be the Perron vector of $A(G')$ and $\rho'=\rho(G')$.
By symmetry, $x$ takes the same value on the vertices of $V(K_{s})$, $V(K_{n-ks-s-1})$ and $V((ks+1)K_{1})$, respectively.
We denote the entry of $x$ by $x_1$, $x_2$ and $x_3$ corresponding to the vertices in the above three vertex sets, respectively.
By $A(G')x=\rho' x$, we have $$\rho'x_3=sx_1.$$ Since $\rho'>0$, then
\begin{equation}\label{for5}
x_3=\frac{sx_1}{\rho'}.
\end{equation}

Let $G^*=K_{\delta}\vee (K_{n-k\delta-\delta-1}+(k\delta+1)K_{1})$.
Suppose that $E_1=\{v_{i}w_{j}|1\leq i\leq n-ks-s-1, k\delta+2\leq j\leq ks+1\}\cup \{w_{i}w_{j}|k\delta+2\leq i\leq ks, i+1\leq j\leq ks+1\}$ and $E_2=\{u_{i}w_{j}|\delta+1\leq i\leq s, 1\leq j\leq k\delta+1\}$. Obviously, $G^*\cong G'+E_1-E_2.$
Note that $A(G^*)$ can be obtained by replacing $s$ with $\delta$ in $A(G')$.
Similarly, let $y$ be the Perron vector of $A(G^*)$ and $\rho^*=\rho(G^*)$.
By symmetry, $y$ takes the same value on the vertices of $V(K_{\delta})$, $V(K_{n-k\delta-\delta-1})$ and $V((k\delta+1)K_1)$, respectively.
We denote the entry of $y$ by $y_1$, $y_2$ and $y_3$ corresponding to the vertices in these three vertex sets, respectively.
By $A(G^*)y=\rho^*y$, we have
\begin{gather}
\rho^*y_2=\delta y_1+(n-k\delta-\delta-2)y_2,\label{for6}\\
\rho^*y_3=\delta y_1.\label{for7}
\end{gather}
Note that $G^*$ contains $K_{n-k\delta-\delta-1}$ as a proper subgraph.
By Lemma \ref{le2}, we have $\rho^*>\rho(K_{n-k\delta-\delta-1})=n-k\delta-\delta-2$.
Putting (\ref{for7}) into (\ref{for6}), and combining $\rho^*>n-k\delta-\delta-2$, we have
\begin{equation}\label{for8}
y_2=\frac{\rho^*y_3}{\rho^*-(n-k\delta-\delta-2)}.
\end{equation}
Recall that $G'=K_{s}\vee (K_{n-ks-s-1}+(ks+1)K_{1})$, which is not a complete graph, so $\rho'<n-1.$
Note that $n\geq ks+s+1$. Then $\delta+1 \leq s \leq \frac{n-1}{k+1}$.

\begin{claim}\label{claim2}
$\rho^*>\rho'$.
\end{claim}

Assume that $\rho'\geq \rho^*$. Since $n\geq ks+s+1$, $k\geq 2$, $x>0$ and $y>0$, we have
\begin{eqnarray*}
&&y^{T}(\rho^*-\rho')x\\
&=&y^{T}(A(G^*)-A(G'))x\\
&=&\sum_{i=1}^{n-ks-s-1}\sum_{j=k\delta+2}^{ks+1}(x_{v_{i}}y_{w_j}+x_{w_j}y_{v_i})
+\sum_{i=k\delta+2}^{ks}\sum_{j=i+1}^{ks+1}(x_{w_i}y_{w_j}+x_{w_j}y_{w_i})
-\sum_{i=\delta+1}^{s}\sum_{j=1}^{k\delta+1}(x_{u_i}y_{w_j}+x_{w_j}y_{u_i})\\
&=&(s-\delta)[k(n-ks-s-1)(x_2y_2+x_3y_2)+k(ks-k\delta-1)x_3y_2-(k\delta+1)(x_1y_3+x_3y_2)]\\
&=&(s-\delta)[(k(n-k\delta-\delta-s-2)-1)x_3y_2+k(n-ks-s-1)x_2y_2-(k\delta+1)x_1y_3]\\
&\geq&(s-\delta)[(k(n-k\delta-\delta-s-2)-1)x_3y_2-(k\delta+1)x_1y_3].
\end{eqnarray*}
By (\ref{for5}) and (\ref{for8}), we obtain that
\begin{eqnarray*}
&&y^{T}(\rho^*-\rho')x\\
&\geq&(s-\delta)x_1y_3\left[\frac{s\rho^*(k(n-k\delta-\delta-s-2)-1)}{\rho'(\rho^*-(n-k\delta-\delta-2))}-(k\delta+1)\right]\\
&=&\frac{(s-\delta)(k\delta+1)x_1y_3}{\rho'(\rho^*-(n-k\delta-\delta-2))}
\left[\frac{s\rho^*(k(n-k\delta-\delta-s-2)-1)}{k\delta+1}-\rho'(\rho^*-(n-k\delta-\delta-2))\right]\\
&=&\frac{(s-\delta)(k\delta+1) x_1y_3}{\rho'(\rho^*-(n-k\delta-\delta-2))}
\left[\frac{s\rho^*(n-k\delta-\delta-s-2-\frac{1}{k})}{\delta+\frac{1}{k}}-\rho'\rho^*+\rho'(n-k\delta-\delta-2)\right]\\
&=&\frac{\rho^*(s-\delta)(k\delta+1)x_1y_3}{\rho'(\rho^*-(n-k\delta-\delta-2))}
\left[\frac{s(n-k\delta-\delta-s-2-\frac{1}{k})}{\delta+\frac{1}{k}}-\rho'+\frac{\rho'}{\rho^*}(n-k\delta-\delta-2)\right]
\end{eqnarray*}
Since $s\geq \delta+1$, $k\geq 2$, $\rho^*\leq \rho'$ and $n\geq3(k+1)\delta+5$, then
\begin{eqnarray*}
&&y^{T}(\rho^*-\rho')x\\
&\geq&\frac{\rho^*(s-\delta)(k\delta+1)x_1y_3}{\rho'(\rho^*-(n-k\delta-\delta-2)}
[(n-k\delta-\delta-s-2-\frac{1}{k})-\rho'+(n-k\delta-\delta-2)]\\
&=&\frac{\rho^*(s-\delta)(k\delta+1)x_1y_3}{\rho'(\rho^*-(n-k\delta-\delta-2))}[2n-2(k+1)\delta-\rho'-s-\frac{1}{k}-4]\\
&\triangleq&\frac{\rho^*(s-\delta)(k\delta+1)x_1y_3}{\rho'(\rho^*-(n-k\delta-\delta-2))}f(n),
\end{eqnarray*}
Note that $s\leq\frac{n-1}{k+1}$, $k\geq2$, $n\geq3(k+1)\delta+5$ and $\rho'<n-1.$ Then
\begin{eqnarray*}
f(n)&=&2n-2(k+1)\delta-\rho'-s-\frac{1}{k}-4\\
&\geq&2n-2(k+1)\delta-\rho'-\frac{n-1}{k+1}-\frac{1}{k}-4\\
&\geq&2n-2(k+1)\delta-\rho'-\frac{n-1}{3}-\frac{1}{2}-4\\
&=&\frac{5n}{3}-2(k+1)\delta-\frac{25}{6}-\rho'\\
&\geq&n+\frac{2}{3}(3(k+1)\delta+5)-2(k+1)\delta-\frac{25}{6}-\rho'\\
&=&n-\rho'-\frac{5}{6}>\frac{1}{6}.
\end{eqnarray*}
Hence $f(n)>0$. Combining $s\geq \delta+1$ and $\rho^*>n-k\delta-\delta-2$, we obtain that $y^{T}(\rho^*-\rho')x>0$.
It follows that $\rho^*>\rho',$ which contradicts the assumption $\rho'\geq\rho^*.$
\hspace*{\fill}$\Box$

By Claim \ref{claim2} and (\ref{for1}), we have
$$\rho(G)\leq\rho(G')<\rho(G^*),$$
a contradiction.
\hspace*{\fill}$\Box$

%\vspace{5mm}
\noindent
%{\bf Declaration of competing interest}
\vspace{3mm}

%The authors declare that they have no known competing financial interests or personal relationships that could have appeared to influence the work reported in this paper.

%\noindent
%{\bf Acknowledgement}

%The authors would like to thank the anonymous referees for their helpful comments on improving the presentation of
%paper.

\end{document}